\documentclass[11pt, reqno]{amsart}

\usepackage{latexsym}
\usepackage{amssymb}
\usepackage{mathrsfs}
\usepackage{amsmath}
\usepackage{fancybox,color}
\usepackage{enumerate}
\usepackage{enumitem}
\usepackage[latin1]{inputenc}
\usepackage{comment}
\usepackage{url}

\newcommand{\kom}[1]{}
\renewcommand{\kom}[1]{{\bf [#1]}}

\renewcommand{\theequation}{\arabic{section}.\arabic{equation}}

 \def\1{\raisebox{2pt}{\rm{$\chi$}}}

\newtheorem{theorem}{Theorem}[section]

\newtheorem{lemma}[theorem]{Lemma}

\newtheorem{definition}[theorem]{Definition}

 \newcommand{\eps}{{\varepsilon}}
 \def\1{\raisebox{2pt}{\rm{$\chi$}}}
 
\newcommand{\abs}[1]{\left|#1\right|}

\def\vint_#1{\mathchoice%
          {\mathop{\kern 0.2em\vrule width 0.6em height 0.69678ex depth -0.58065ex
                  \kern -0.8em \intop}\nolimits_{\kern -0.4em#1}}%
          {\mathop{\kern 0.1em\vrule width 0.5em height 0.69678ex depth -0.60387ex
                  \kern -0.6em \intop}\nolimits_{#1}}%
          {\mathop{\kern 0.1em\vrule width 0.5em height 0.69678ex
              depth -0.60387ex
                  \kern -0.6em \intop}\nolimits_{#1}}%
          {\mathop{\kern 0.1em\vrule width 0.5em height 0.69678ex depth -0.60387ex
                  \kern -0.6em \intop}\nolimits_{#1}}}
\def\vintslides_#1{\mathchoice%
          {\mathop{\kern 0.1em\vrule width 0.5em height 0.697ex depth -0.581ex
                  \kern -0.6em \intop}\nolimits_{\kern -0.4em#1}}%
          {\mathop{\kern 0.1em\vrule width 0.3em height 0.697ex depth -0.604ex
                  \kern -0.4em \intop}\nolimits_{#1}}%
          {\mathop{\kern 0.1em\vrule width 0.3em height 0.697ex depth -0.604ex
                  \kern -0.4em \intop}\nolimits_{#1}}%
          {\mathop{\kern 0.1em\vrule width 0.3em height 0.697ex depth -0.604ex
                  \kern -0.4em \intop}\nolimits_{#1}}}

\newcommand{\aveint}[2]{\mathchoice%
          {\mathop{\kern 0.2em\vrule width 0.6em height 0.69678ex depth -0.58065ex
                  \kern -0.8em \intop}\nolimits_{\kern -0.45em#1}^{#2}}%
          {\mathop{\kern 0.1em\vrule width 0.5em height 0.69678ex depth -0.60387ex
                  \kern -0.6em \intop}\nolimits_{#1}^{#2}}%
          {\mathop{\kern 0.1em\vrule width 0.5em height 0.69678ex depth -0.60387ex
                  \kern -0.6em \intop}\nolimits_{#1}^{#2}}%
          {\mathop{\kern 0.1em\vrule width 0.5em height 0.69678ex depth -0.60387ex
                  \kern -0.6em \intop}\nolimits_{#1}^{#2}}}
\newcommand{\ud}{\, d}
\newcommand{\half}{{\frac{1}{2}}}

\newcommand{\Om}{\Omega}

\newcommand{\om}{\omega}

\newcommand{\diam}{\operatorname{diam}}

\newcommand{\ber}{\operatorname{Ber}}

\begin{document}

\title[Uniform measure density condition and game regularity]{Uniform measure density condition and game regularity for tug-of-war games}

\author[Joonas Heino]{Joonas Heino}
\address{Department of Mathematics and Statistics, University of Jyvaskyla, PO~Box~35, FI-40014 Jyvaskyla, Finland}
\email{joonas.p.heino@student.jyu.fi}

\date{\today}

\begin{abstract}
We show that a uniform measure density condition implies game regularity for all $2< p<\infty$ in a stochastic game called 'tug-of-war with noise'. The proof utilizes suitable choices of strategies combined with estimates for the associated stopping times and density estimates for the sum of independent and identically distributed random vectors. 
\end{abstract}

\maketitle

\section{Introduction}
The profound connection between stochastic processes and classical linear partial differential equations has been pivotal. For example, this connection was made use of in \cite{krylovs79,krylovs80} to establish regularity result for the second order equations in a non divergence form. Recently, a connection between nonlinear infinity harmonic functions and tug-of-war games was discovered in \cite{peresssw09}. Later in \cite{peress08}, the authors found a stochastic game related to $p$-harmonic functions. They proved among other things by using a game approach that in a game regular domain there exists a $p$-harmonic function extending continuously to the boundary with the given continuous boundary values. However, a problem asking if a regular boundary point for the $p$-Laplacian is necessarily game regular was left open.

We study a modified version of a 'tug-of-war with noise' developed in \cite{manfredipr12} and also related to $p$-harmonic functions. First, the players choose a step length $\epsilon>0$ and a starting point $x_0$. Then, they toss a biased coin, and if they get heads (probability $\alpha$), the players play a 'tug-of-war', that is, they toss a fair coin and the winner of the toss can move the game position to any point of the open ball centered at $x_0$ and of the radius $\epsilon$. If in the first toss, they get tails (probability $\beta$), the game point moves according to the uniform distribution in the open ball centered at $x_0$ and of the radius $\epsilon$. After the first move, the players play the same game from the new game position. The game ends, when the game position exits the game domain for the first time. In the end, Player 2 pays to Player 1 the amount given by the payoff function at the first point outside the domain. We consider this version of the game because the players do not affect the direction of the noise and hence, we can prove sharp enough estimates for the density of the noise. 

We give a stochastic proof that a uniform measure density condition implies game regularity (Theorem \ref{paatulos p suurempi kaksi}). Roughly, a boundary point $y$ is game regular, if Player 1 has a strategy to end the game near $y$ with a probability close to one whenever the game starts near $y$ as well. A boundary point $y$ satisfies a measure density condition, if the Lebesgue measure of the complement of the game domain in the ball centered at $y$ is comparable to the Lebesgue measure of the whole ball. The proof of Theorem \ref{paatulos p suurempi kaksi} utilizes a stochastic density estimate for the sum of independent and identically distributed random vectors (Lemma \ref{ehka kiinnostava}). In addition, we use a 'cylinder walk' framework together with a cancellation strategy for Player 1 to connect the stochastic estimates to the setting. We omit the case $p=2$, because in that case the process is merely a random walk and the result follows from the classical invariance principle. 

Game theory has already given new insights to partial differential equations. For instance, the ideas emerging from nonlinear game theory have led to simpler as well as completely different proofs for PDEs (see for example \cite{scottc10} and \cite{luirops13}). In addition, a dynamic programming principle related to the game also arises from discretization schemes (see for instance \cite{oberman05}).

We expect the techniques developed in this paper to be useful for a larger class of partial differential equations as well. In addition, stochastic estimates on where the game spends time under cancellation strategies are likely to be important for further results.

This work is organized as follows. In Section 2, we describe the preliminaries needed in the paper. Then in Section 3, we show that the uniform measure density condition implies game regularity for all $2< p<\infty$. For brevity, we do not write down all the stochastic calculations needed in the section, but the calculations are in the appendix.

\section{Preliminaries}
First, let us start by introducing the notation. We denote the standard Euclidian open ball by $B_r(x_0)\subset \mathbb R^n$,
$$
B_r(x_0)=\{z \in \mathbb R^n: |z-x_0|<r \}.
$$
Lebesgue measure is denoted by $|\cdot|$, and in addition, the notation $C_{n,p}$ means that the universal constant depends only on $n$ and $p$. Throughout the paper, we use the asymptotic notation $\mathcal O(\epsilon)$. For example, if a real-valued function $f$ satisfies the inequality $f(\epsilon)\leq \mathcal O (\epsilon)$, it means that there exists a constant $C>0$ such that $|f(\epsilon)|\leq C \epsilon$ for all $\epsilon>0$ small enough. 

Let $2< p <\infty$, $\epsilon>0$ and dimension $n\geq 1$. Fix a bounded, non-empty and open set $\Om\subset \mathbb R^n$. Next, we recall the two-player zero-sum-game called 'tug-of-war with noise'. First, choose a starting point $x_0\in \Omega$ for the game, and then, the players toss a biased coin with probabilities $\alpha$ and $\beta$. The probabilities depend on $n$ and $p$ by
\begin{align}\label{alpha ja beeta}
\alpha = \frac{p-2}{p+n},~~~~~\beta=\frac{n+2}{p+n}.
\end{align}
The players get heads with the probability $\alpha$, and in this case, they will toss a fair coin and the winner of the toss can move the game position to any point of the open ball $B_\epsilon(x_0)$. Tossing of a fair coin and the movement after the toss are the 'tug-of-war' parts of the game. On the other hand, if they get tails, the next game position will be decided by the uniform distribution in the ball $B_\epsilon(x_0)$. A random movement is the 'noise' part of the game. After the first move is decided, the players continue playing the same game from the new position.

The game procedure yields a sequence of game positions $x_0,x_1,x_2,\dots$, where every $x_k$ is a random variable. A history of a game up to step $k$ is a vector of the first $k+1$ game positions $x_0,\dots,x_k$ and $k$ coin tosses $c_1,\dots,c_k$, that is,
$$
h_k:=(x_0,(c_1,x_1),\dots,(c_k,x_k)).
$$
In the above, $c_j\in \mathcal C:=\{0,1,2\}$, where 0 denotes that Player 1 wins, 1 that Player 2 wins and 2 that a random movement occurs.

To prescribe boundary values, let us denote a compact boundary strip of width $\epsilon$ by
$$
\Gamma_\epsilon:=\{z\in \mathbb R^n\setminus \Omega: \inf_{y\in\partial \Omega}|z-y|\leq \epsilon\}.
$$ 
The reason to use the boundary strip instead of just the boundary is that $B_\epsilon(x) \subset \Omega_\epsilon:=\Omega \cup \Gamma_\epsilon$ for all $x\in \Omega$. After the first time the game position is in $\Gamma_\epsilon$, the players do not move it anymore. For all $k\geq 0$, the history $h_k$ belongs to the space $H^k:=x_0\times (\mathcal C,\Omega_\epsilon)^k$ with $H^0:=x_0$. We denote the space of all game sequences by
$$
H^\infty:=\bigcup_{k\geq 0} H^k=x_0\times (\mathcal C,\Omega_\epsilon) \times (\mathcal C,\Omega_\epsilon) \times \cdots.
$$

A \textit{strategy} for Player 1 is a sequence of Borel measurable functions that give the next game position given the history of the game. To be more precise, a strategy for Player 1 is $S_1:=(S_{1,k})_{k=0}^\infty$ with
$$
S_{1,k}:H^k\to \mathbb R^n
$$
for all $k\geq0$. For example, if Player 1 wins the $(k+1)$th toss,
$$
S_{1,k}(x_0,(c_1,x_1),\dots,(c_k,x_k))=x_{k+1}\in B_\epsilon(x_k)
$$
for all $h_k \in H^k$. Similarly Player 2 deploys a strategy $S_2$.

We denote the first hitting time to the set $\Gamma_\epsilon$ by
$$
\tau:=\tau(\omega)=\inf\{k: x_k \in \Gamma_\epsilon,k=0,1,2,\dots\}.
$$
The game process is a discrete time adapted process with respect to the filtration $\mathcal F_0:=\sigma(x_0)$ and
$$
\mathcal F_k:=\sigma(x_0,(c_1,x_1),\dots,(c_k,x_k))~\text{for } k\geq1,
$$
so $\tau$ is always a stopping time. The game ends at the random time $\tau$, and the payoff is $F(x_\tau$), where $F:\Gamma_\epsilon \to \mathbb R$ is a fixed, bounded and Borel measurable \textit{payoff function}. In the end, Player 2 pays the amount $F(x_\tau)$ to Player 1.

To establish a unique probability measure, we need to know a starting point $x_0$ and strategies $S_1$ and $S_2$. Then, the probability measure $\mathbb P_{S_1,S_2}^{x_0}$ on the natural product $\sigma$-algebra is built by applying Kolmogorov's extension theorem to the family of transition densities
\begin{align*}
&\pi_{S_1,S_2}(x_0,(c_1,x_1),\dots,(c_k,x_k),(C,A)) \\
&=\frac{\alpha}{2}\delta_0(C)\delta_{S_1(x_0,(c_1,x_1),\dots,(c_k,x_k))}(A)+\frac{\alpha}{2}\delta_1(C)\delta_{S_2(x_0,(c_1,x_1),\dots,(c_k,x_k))}(A) \\
&\text{  }+\beta \delta_2(C)\frac{|A\cap B_\epsilon(x_k)|}{|B_\epsilon(x_k)|}
\end{align*}
for any subset $C\subset \mathcal C$ and Borel subset $A\subset \Omega_\epsilon$ as long as $x_k \in \Omega$. If $x_k \not\in \Omega$, the transition probability forces $x_{k+1}=x_k$.

The expected payoff is 
$$
\mathbb E_{S_1,S_2}^{x_0}[F(x_\tau)]=\int_{H^\infty} F(x_\tau(\omega)) \ud \mathbb P_{S_1,S_2}^{x_0},
$$
when the game starts from $x_0$ and the players use strategies $S_1$ and $S_2$. The \textit{value of the game for Player 1} is given by
$$
u_\epsilon^1(x_0)=\sup_{S_1} \inf_{S_2} \mathbb E_{S_1,S_2}^{x_0}[F(x_\tau)]
$$ 
and the \textit{value of the game for Player 2} is given by
$$
u_\epsilon^2(x_0)=\inf_{S_2} \sup_{S_1} \mathbb E_{S_1,S_2}^{x_0}[F(x_\tau)],
$$
respectively. The game has a value i.e. there exists a unique value function $u_\epsilon:=u_\epsilon^1=u_\epsilon^2$ (see \cite{manfredipr12} and \cite{luirops14}).

Since $\Omega$ is bounded, the game ends almost surely for any choice of strategies. This is true due to the fact that for $n_0\geq 1 $ large enough, we have $n_0\epsilon>\diam(\Omega)$, and almost surely there will be infinitely many blocks of length $n_0$ consisting of solely random moves in the game. 

Observe that the history $h_k$ contains all the information at the moment $k$, and since the strategies are a collection of Borel measurable functions from all possible histories, it is clear that the game process will not be a \textit{Markov process} in general. 

This version of the tug-of-war game has good symmetry properties, which we will utilize in the proofs. Other versions of tug-of-war games have been studied for example in \cite{peress08} and \cite{kawohlmp12} and a continuous time game in \cite{atarb10}. 

A rough outline of the connection between the version of the game considered in this paper and $p$-harmonic functions is the following. First, assume that we have a $p$-harmonic function in an open set $\Omega'\supset \Omega$ with a nonvanishing gradient. Then, the $p$-harmonic function is real analytic, and Theorem 4.1 in \cite{manfredipr12} states that the game with probabilities \eqref{alpha ja beeta} and with the values of the $p$-harmonic function on the boundary approximates the $p$-harmonic function in the game domain. The proof is based on the gradient strategy for the $p$-harmonic function and on the optional stopping theorem as well as on the asymptotic expansion in \cite{manfredipr10}. 

The general case requires game regularity of the boundary of the game domain. Then, it is possible to use a barrier argument to get estimates close to the boundary. By copying the strategies and utilizing the translation invariance of the game, the same estimates also holds in the interior of the game domain. Finally, a variant of the classical Arzel\`{a}-Ascoli's theorem provides a convergent subsequence. To prove that the limit is a viscosity solution to the homogeneous $p$-Laplace equation, a dynamic programming principle related to the game is applied (for more details about the principle, see for example \cite{luirops14}).

\section{Measure density condition implies game regularity}

We show in Theorem \ref{paatulos p suurempi kaksi} that a uniform measure density condition implies game regularity for all $p>2$. To establish this, we first show in Lemma \ref{ehto} a more attainable criterion for game regularity. Then in Theorem \ref{sylinteri arvio}, we use a 'cylinder walk' framework, introduced in \cite{luirops13}, to obtain some important hitting probability estimates.

\begin{definition}

A point  $y\in \partial \Om$ satisfies a measure density condition if there is $c>0$ such that
\[
\begin{split}\label{eq:measure-density}
\abs{\Om^c\cap B_r(y)}\ge c \abs{B_r(y)}
\end{split}
\]
for all $r>0$.
\end{definition}

\begin{definition}\label{game-regular}
A point $y\in \partial \Omega$ is game regular, if for all $\delta>0$ and $\eta>0$, there exist $\delta_0>0$ and $\epsilon_0>0$ such that for all fixed $\epsilon<\epsilon_0$ and for all $x_0 \in B_{\delta_0}(y)$, there is a strategy $S^*_1$ for Player 1 such that
$$
\mathbb P_{S_1^*,S_2}^{x_0}(x_\tau \in B_\delta(y)\cap \Omega^c)\geq 1-\eta. 
$$ 
If every boundary point of $\Omega$ is game regular, we say that $\Omega$ is game regular.
\end{definition}

Roughly speaking, game regularity means that whenever the game starts near a boundary point $y$, Player 1 has a strategy to end the game near $y$ with a high probability. Next, we give a more attainable criterion to obtain game regularity. We modify the idea from \cite[p.\,13]{peress08}.

\begin{lemma}\label{ehto}
 A boundary point $y\in \partial \Omega$ is game regular if there exists a constant $\theta>0$ such that for all $\delta>0$, there are parameters $\epsilon_0>0$ and $\delta_0>0$ such that for all fixed $\epsilon<\epsilon_0$ and for all $x_0\in B_{\delta_0}(y)$, there is a strategy $S_1^*$ for Player 1 such that
$$
\mathbb P_{S_1^*,S_2}^{x_0}(\text{the game ends before exiting the ball }B_\delta(y))\geq \theta.
$$

\begin{proof}
The idea of the proof is the following. By choosing $\delta_0>0$ small enough, we can start the game as near the point $y$ as we want, and in order to exit the ball $B_\delta(y)$, the game sequence has to exit all the concentric smaller balls inside $B_\delta(y)$ as well. The probability to exit all the concentric balls inside $B_\delta(y)$ can be estimated above via the uniform probability $\theta$; it is less than $(1-\theta)^k$, where $k$ is the amount of concentric balls inside $B_\delta(y)$. Thus, the probability to end the game near $y$ is close to one, when $k$ is big enough.

To be more precise, let $\delta>0$ and $\eta>0$. Now, there are $\theta>0$, $\epsilon_{0,1}>0$ and  $0<\delta_{0,1}<\delta$ such that for all $\epsilon<\epsilon_{0,1}$ and for all $x_0\in B_{\delta_{0,1}}(y)$, we have a strategy $S_1^1$ for Player 1  such that
$$
\mathbb P_{S_1^1,S_2}^{x_0}(\text{the game ends before exiting the ball }B_\delta (y))\geq \theta.
$$
We can assume that $\epsilon_{0,1}<\delta_{0,1}$. Again similarly as above, for the constant $\delta_{0,1}-\epsilon_{0,1}$, there are $\epsilon_{0,2}>0$ and $0<\delta_{0,2}<\delta_{0,1}/2$ such that for all $\epsilon<\epsilon_{0,2}$ and for all $x_0\in B_{\delta_{0,2}}(y)$, we have a strategy $S_1^2$ for Player 1 such that the probability to end the game before exiting the ball $B_{\delta_{0,1}-\epsilon_{0,1}}(y)$ is at least $\theta$. We can do this as many times we want. Let us do this $k\in \mathbb N$ times, where $k$ is such that  
$$
(1-\theta)^k \leq \eta.
$$ 

Define $\delta_0:=\delta_{0,k}$ and $\epsilon_0:=\min \{\epsilon_{0,1},\dots,\epsilon_{0,k}\}$, and fix any $x_0 \in B_{\delta_0}(y)$ and $\epsilon<\epsilon_0$. We can assume that $\epsilon<\min\{\delta_0, \delta_{0,k-1}-\delta_{0,k},\dots,\delta_{0,1}-\delta_{0,2}\}$ so that the game position cannot jump over many concentric balls during one turn. Denote the first time the game sequence exits $B_{\delta_{0,i-1}-\epsilon_{0,i-1}}(y)$ by $\tau^i:=\tau^i(\omega)$ for all $i\in\{1,\dots,k\}$ with $\delta_{0,0}:=\delta$ and $\epsilon_{0,0}:=0$. Also, denote the set
\begin{align*}
A_i&:= \{\text{exits the ball $B_{\delta_{0,i-1}-\epsilon_{0,i-1}}(y)$ before the game ends}\} 
\end{align*}
for all $i\in \{1,\dots k\}$.

Recall that the game ends at the random time $\tau$. Define a strategy $S_1^*$ for Player 1 such that first, Player 1 uses the strategy $S_1^k$. If $\tau^k<\tau$, Player 1 starts to use the strategy $S_1^{k-1}$ after the stopping time $\tau^k$. Similarly, if $\tau^{k-1}<\tau$, Player 1 starts to use the strategy $S_1^{k-2}$ after the stopping time $\tau^{k-1}$. Thus, if we have stopping times $0<\tau^k< \tau^{k-1}< \dots < \tau^1<\tau$, after every stopping time $\tau^i$, Player 1 starts to use the strategy $S_1^{i-1}$ for all $i\in\{2,\dots,k\}$ and for all game sequences $\omega\in H^\infty$. After the stopping time $\tau^1$, Player 1 does not change her strategy anymore. Observe that the earlier strategy $S_1^{i}$ does not affect the game after the first time the game sequence exits $B_{\delta_{0,i-1}-\epsilon_{0,i-1}}(y)$ for every $i\in\{2,\dots,k\}$. Roughly this means that for every $i\in\{2,\dots,k\}$, after the stopping time $\tau^i$, Player 1 forgets everything that has happened prior the time $\tau^i$.

Let $S_2$ be any strategy for Player 2. The strategy $S_2$ can depend heavily on the past, so it could well be that our game process does not have any Markovian structure at any game round. However, the uniform $\theta$ is independent of the information available, so roughly, Player 2 cannot gain too much from the information of the past.

By the reasoning above, we can estimate iteratively
\begin{align*}
&\mathbb P^{x_0}_{S_1^*,S_2}(\text{exits the ball $B_\delta(y)$ before the game ends}) \\
=&\mathbb E^{x_0}_{S_1^k,S_2}\bigg[\1_{A_k}\mathbb E^{x_0}_{S_1^{k-1},S_2}\bigg[\prod_{l=1}^{k-1}\1_{A_l}\Big|\mathcal F_{\tau^k}\bigg]\bigg] \\
= &  \mathbb E^{x_0}_{S_1^k,S_2}\bigg[\1_{A_k}\mathbb E^{x_0}_{S_1^{k-1},S_2}\bigg[\1_{A_{k-1}}\cdots \mathbb E^{x_0}_{S_1^{1},S_2}\Big[\1_{A_1}\big|\mathcal F_{\tau^2}\Big]\cdots \Big|\mathcal F_{\tau^{k}}\bigg]\bigg] \\
\leq & (1-\theta)^k \\
\leq & \eta.
\end{align*}
This implies that
\begin{align*}
\mathbb P_{S_1^*,S_2}^{x_0}(x_\tau \in B_\delta(y)\cap \Omega^c) &\geq \mathbb P_{S_1^*,S_2}^{x_0}(\text{the game ends before exiting }B_\delta(y))\\
&\geq 1-\eta
\end{align*}
i.e. we have shown the game regularity.

\end{proof} 
\end{lemma}

To see that the uniform measure density condition implies game regularity, we need a 'cylinder walk' framework. 

\textbf{Cylinder walk}. Set the constants $\alpha$,$\beta>0$ with $\alpha+ \beta=1$ as before in \eqref{alpha ja beeta}, and fix the cylinder size $r>0$. Consider the following random walk (called the 'cylinder walk') in a $n+1$ -dimensional cylinder $B_{r}(0)\times [0,r]$. Suppose that we are at a point $(x_j,t_j)\in B_r(0) \times [0,r]$. Next, we move to the point $(x_j,t_j-\epsilon)$ with probability $\alpha/2$ and to the point $(x_j,t_j+\epsilon)$ with probability $\alpha /2$. With probability $\beta$ we move to the point $(x_{j+1},t_j)$, where $x_{j+1}$ is chosen from the ball $B_\epsilon(x_j)$ according to the uniform distribution. 

We have the following estimate for the probability that the cylinder walk exits the cylinder through its bottom; the proof is in the appendix of the paper \cite{luirops13}. 

\begin{lemma}\label{cylinder}
Let us start the cylinder walk from the point $(0,t)$ with $0<t<r$. Then, the probability that the walk exits the cylinder through its bottom is at least
$$
1-C_{n,p}(t+\epsilon)/r
$$
for all $\epsilon>0$ small enough.
\end{lemma}

Assume that the origin $0 \in \mathbb R^{n+1}$ at the bottom of the cylinder belongs to the set $\partial \Omega \times \{0\}$ and that this boundary point satisfies the measure density condition. The set $\Omega \cap B_r(0) \times \{0\} \subset B_r(0)\times [0,r]$. We are interested in the probability that the cylinder walk exits through the bottom and in addition, at the first time the walk hits the bottom, the process is in the complement of the set $\Omega$. Since the origin satisfies the measure density condition, the complement has some positive Lebesgue measure. This suggests that the event we are interested in could have some positive probability measure. 

The cylinder walk can be constructed by combining three independent random constructions. There is a 'horizontal' random walk with the initial position $\tilde{x}_0=x\in B_r(0)$. The point $\tilde{x}_{j+1}$ is chosen according to the uniform distribution in the ball $B_\epsilon(\tilde{x}_j)\subset \mathbb R^n$ for all $j\geq 0$. Further, there is a 'vertical' random walk in the real axis with steps $+\epsilon$ or $-\epsilon$ and with the initial position $\tilde{t}_0=t\in ]0,r[$. For all $j\geq 0$, the next positions are $\tilde{t}_{j+1}=\tilde{t}_{j}+\epsilon$ or $\tilde{t}_{j+1}=\tilde{t}_{j}-\epsilon$ both with probability $\frac{1}{2}$. In addition, there is the increasing sequence
$$
U_j=\sum_{m=1}^j \ber_m,
$$
where the $\ber_m$:s are independent Bernoulli variables with $\ber_m(\omega) \in \{0,1\}$ and $\mathbb P(\ber_m=1)=\alpha$. Therefore, a copy of the cylinder walk is obtained by letting for $j\geq 0$
$$
t_j=\tilde{t}_{U_j},~~~x_j=\tilde{x}_{j-U_j}.
$$

Let $\tau_g$ stand for the first moment $t_j$ exits the cylinder through its bottom or top i.e. the first $j$ such that $t_j\in \mathbb R \setminus  ]0,r[$. Also, let $\tilde{\tau}_g$ stand for the first moment $\tilde{t}_j$ exits the cylinder through its bottom or top. Here, the subindex $g$ refers to a 'good exit'.

We assume that the walk starts from the point $(0,t)$ with $0<t<r$ i.e. $\tilde{x}_0=0$ and $\tilde{t}_0=t$. First, let us study the properties of the function $\tau_g-U_{\tau_g}=\tau_g-\tilde{\tau}_g$. The random variable $\tau_g-\tilde{\tau}_g$ is the number of times a random horizontal movement has occured at the first moment the cylinder walk hits the bottom. The proof of the lemma below is in the appendix for completeness.

\begin{lemma}\label{LEMMA tod nak}
Let $\tau_g,\tilde{\tau}_g,\alpha,\beta, t$ and $r$ be as above. The random variable $\tau_g-\tilde{\tau}_g$ holds the following properties for all $a>0$
\begin{align}
&\mathbb P(\tau_g-\tilde{\tau}_g \geq a\epsilon^{-1}) \geq 1 -\mathcal O (\epsilon) \text{ and} \label{toka tulos kappale}\\
&\mathbb P(\tau_g-\tilde{\tau}_g \geq a\epsilon^{-2})\leq  1- \frac{2}{\sqrt{2 \pi}} \int_{\frac{\min\{t,r-t\}}{\sqrt{a}}\nu_{n,p}}^\infty e^{-\frac{s^2}{2}}\ud s+\mathcal O(\epsilon)\label{kolmas tulos kappale}
\end{align}
with the constant
$$
\nu_{n,p}:=2\sqrt{\frac{\beta+0.01\alpha}{0.99\alpha}}.
$$
\end{lemma}
For any $a>0$, the inequalities \eqref{toka tulos kappale} and \eqref{kolmas tulos kappale} yield 
\begin{align}
\begin{split}\label{oikea skaala}
\mathbb P (a\epsilon^{-1} \leq \tau_g-\tilde{\tau}_g <a\epsilon^{-2}) &\geq \mathbb P ( \tau_g-\tilde{\tau}_g \geq a\epsilon^{-1} )- \mathbb P(\tau_g-\tilde{\tau}_g \geq a\epsilon^{-2}) \\
& \geq \frac{2}{\sqrt{2\pi}}\int_{\frac{\min\{t,r-t\}}{\sqrt{a}}\nu_{n,p}}^\infty e^{-\frac{s^2}{2}} \ud s-\mathcal O(\epsilon).
\end{split}
\end{align} 
Observe that 
$$
\frac{2}{\sqrt{2\pi}}\int_{\frac{\min\{t,r-t\}}{\sqrt{a}}\nu_{n,p}}^\infty e^{-\frac{s^2}{2}} \ud s \to 1
$$
as $t \to 0$. Thus, the inequality \eqref{oikea skaala} points out that for the cylinder walk started near $(0,0)$, the random variable $\tau_g-\tilde{\tau}_g$ is very likely between the times $a\epsilon^{-1}$ and $a\epsilon^{-2}$ for any $\epsilon$ small enough and fixed $a>0$.

Next, we concentrate on the distribution of the random variable $\tilde{x}_k$. Assume that $Z$ is a random vector with the uniform distribution in the ball $B_\epsilon(0)\subset \mathbb R^n$. The density of the random vector $Z$ is 
$$
f_Z(x) = \frac{1}{|B_\epsilon(0)|} \1_{B_\epsilon(0)}(x).
$$
We denote the measure of the unit ball by $\omega_n:=|B_1(0)|$. Let $k_0:=k_{0,n}>2$ denote the constant in Lemma \ref{ehka kiinnostava} and fix any $k\geq k_0$. For the density of the random variable $\tilde{x}_k=\sum_{i=1}^k Z_i$, where the random vectors $Z_i$ are independent and distributed as $Z$, we use the notation $f_k:=f_{\sum_{i=1}^k Z_i}$. The density $f_k$ is a decreasing radial function. In the appendix, we have derived in \eqref{lopullinen arvio origo} and \eqref{lopullinen} the following estimates: There are constants $C_n>0$ and $C_1>0$ such that 
\begin{align*}
f_k(0) \leq C_n \bigg(\frac{1}{\sqrt{k}\epsilon}\bigg)^n,
\end{align*}
and 
\begin{equation}\label{lopullinen kappale}
f_k\big(C_*\sqrt{k}\epsilon\big) \geq \bigg(\frac{1}{C_1}\bigg)^n\bigg(\frac{0.99}{\omega_n}-C_n(C_*)^n\bigg)\bigg(\frac{1}{\sqrt{k}\epsilon}\bigg)^n
\end{equation}
for all $C_*\in]0,C_1[$. By the comment after the statement of Lemma \ref{ehka kiinnostava} in the appendix, we have
$$
f_k\big(C_*\sqrt{k}\epsilon\big)\geq \zeta\bigg(\frac{1}{\sqrt{k}\epsilon}\bigg)^n
$$
for some $\zeta:=\zeta_n>0$, if we choose $C_*>0$ so small that
\begin{equation}\label{vakion valinta kappale}
C_*<\bigg(\frac{0.99}{\omega_nC_n}\bigg)^{1/n}.
\end{equation}

Let $\tau_b$ stand for the first $j$ when $|x_j|$ reaches $[r,\infty[$. Here, the subindex $b$ refers to a 'bad exit'. Recall that the origin at the bottom of the cylinder satisfies the measure density condition. Let $C_{n,p}>0$ denote the constant in Lemma \ref{cylinder}, and for all $\delta>0$, denote
$$
A_\delta:=B_{\delta}(0) \cap \Omega^c.
$$
\begin{theorem}\label{sylinteri arvio}
Consider the cylinder $B_{\delta/3}(0)\times [0,\delta/3]$ for any fixed $\delta>0$. Then, there exist constants $\theta:=\theta_{n,p}>0,\epsilon_0:=\epsilon_{0,n,p,\delta}>0$ and $\delta_0:=\delta_{0,n,p,\delta}>0$ such that
$$
\mathbb P(\tau_b \leq \tau_g ~\text{or }t_{\tau_g} \geq \delta/3~\text{or }x_{\tau_g} \not \in A_{\delta/3}) \leq 1-\theta
$$
for all $\epsilon<\epsilon_0$ whenever the cylinder walk starts from the point $(0,t)$ for some $0<t\leq \delta_0$.
\begin{proof}
To establish the result, we use the inequality \eqref{oikea skaala} to estimate how many times it is likely that a random horizontal movement has occured at the first time the cylinder walk hits the bottom. Then, we use the estimate \eqref{lopullinen kappale} and the fact that the vertical and horizontal movements are independent to estimate the probability that we are in the complement of the set $\Omega$ at the first time the walk exits the cylinder through its bottom.

Let $0<\lambda<1$, where the exact value of $\lambda$ will be fixed later. Define
\begin{align}\label{delta nolla}
\delta_0:=\frac{\delta\lambda}{3C_{n,p}},
\end{align}
and start the cylinder walk from the point $(0,t)$ for some $0<t\leq \delta_0$ in the cylinder $B_{\delta/3}(0)\times [0,\delta/3]$.

Lemma \ref{cylinder} states that 
$$
\mathbb P( \tau_b \leq \tau_g ~\text{or }t_{\tau_g} \geq \delta/3)\leq 3C_{n,p}\delta^{-1}(t+\epsilon) \leq 3C_{n,p}\delta^{-1}(\delta_0+\epsilon).
$$
Therefore, we have by \eqref{delta nolla} that
$$
\mathbb P(\tau_b \leq \tau_g ~\text{or }t_{\tau_g} \geq \delta/3~\text{or }x_{\tau_g} \not \in A_{\delta/3}) \leq \mathcal O(\epsilon)+\lambda+1-\mathbb P(x_{\tau_g} \in A_{\delta/3}).
$$
The inequality \eqref{oikea skaala} and the remark after suggest the estimate
\begin{align*}
\mathbb P(x_{\tau_g} \in A_{\delta/3})& =\mathbb P(\tilde{x}_{\tau_g-\tilde{\tau}_g} \in A_{\delta/3}) \\
&\geq \mathbb P\Big(\tilde{x}_{\tau_g-\tilde{\tau}_g} \in A_{\delta/3} \text{ and } \delta^2\epsilon^{-1} \leq \tau_g-\tilde{\tau}_g <\delta^2\epsilon^{-2}\Big) \\
&=\sum_{k=\lceil \delta^2\epsilon^{-1}\rceil}^{\lfloor \delta^2\epsilon^{-2}\rfloor} \mathbb P\Big(\tilde{x}_{\tau_g-\tilde{\tau}_g} \in A_{\delta/3} \text{ and } \tau_g-\tilde{\tau}_g =k\Big).
\end{align*}
Denote the index set
$$
I :=\{\lceil \delta^2\epsilon^{-1}\rceil,\lceil \delta^2\epsilon^{-1}\rceil+1,\dots, \lfloor \delta^2\epsilon^{-2}\rfloor \}.
$$
Since the random variables $\tilde{x}_k$ and $\tau_g-\tilde{\tau}_g$ are independent for all $k\in I$, we have
\begin{align*}
&\sum_{k\in I}\mathbb P\Big(\tilde{x}_{\tau_g-\tilde{\tau}_g} \in A_{\delta/3} \text{ and } \tau_g-\tilde{\tau}_g =k\Big)\\
&=\sum_{k\in I} \mathbb P\Big(\tilde{x}_{k} \in A_{\delta/3}\Big)\mathbb P \Big(  \tau_g-\tilde{\tau}_g =k\Big).
\end{align*}
Let $k \in I$ and choose the constant $C_*>0$ as in \eqref{vakion valinta kappale}. We may assume that $C_*<1/3$ and $k\geq k_0$, where $k_0>2$ is the constant in Lemma \ref{ehka kiinnostava}. Because $\sqrt{k}\epsilon <\delta$ and the density $f_k$ is a decreasing radial function, we can calculate
\begin{align*}
\mathbb P\big(\tilde{x}_{k} \in A_{\delta/3}\big)& \geq \mathbb P\big(\tilde{x}_{k} \in B_{C_*\sqrt{k}\epsilon}(0)\cap \Omega^c\big) \\
&\geq f_k\big(C_*\sqrt{k}\epsilon\big)|B_{C_*\sqrt{k}\epsilon}(0) \cap \Omega^c|.
\end{align*}
By using the estimate \eqref{lopullinen kappale} and the uniform measure density condition, we obtain
\begin{align*}
& f_k\big(C_*\sqrt{k}\epsilon\big)|B_{C_*\sqrt{k}\epsilon}(0) \cap \Omega^c| \\
&\geq \bigg(\frac{1}{C_1}\bigg)^n\bigg(\frac{0.99}{\omega_n}-C_n(C_*)^n\bigg)\bigg(\frac{1}{\sqrt{k}\epsilon}\bigg)^nc|B_{C_*\sqrt{k}\epsilon}(0)| \\
&= \omega_nc\bigg(\frac{C_*}{C_1}\bigg)^n\bigg(\frac{0.99}{\omega_n}-C_n(C_*)^n\bigg),
\end{align*}
where the constant $c>0$ comes from the uniform measure density condition. This together with the inequality \eqref{oikea skaala} yield 
\begin{align*}
&\sum_{k\in I}\mathbb P\Big(\tilde{x}_{k} \in A_{\delta/3}\Big)\mathbb P \Big(  \tau_g-\tilde{\tau}_g =k\Big) \\
&\geq   \omega_nc\bigg(\frac{C_*}{C_1}\bigg)^n\bigg(\frac{0.99}{\omega_n}-C_n(C_*)^n\bigg) \mathbb P \Big(  \delta^2\epsilon^{-1} \leq \tau_g-\tilde{\tau}_g <\delta^2\epsilon^{-2}\Big) \\
&\geq \omega_nc\bigg(\frac{C_*}{C_1}\bigg)^n\bigg(\frac{0.99}{\omega_n}-C_n(C_*)^n\bigg)\frac{2}{\sqrt{2\pi}}\int_{\frac{\lambda}{3}\tilde{C}_{n,p}}^\infty e^{-\frac{s^2}{2}}\ud s-\mathcal O(\epsilon)
\end{align*}
with the constant
$$
\tilde{C}_{n,p}:=\frac{2}{C_{n,p}}\sqrt{\frac{\beta+0.01\alpha}{0.99\alpha}}.
$$
Define
$$
\theta_{n,p}:=\omega_nc\bigg(\frac{C_*}{C_1}\bigg)^n\bigg(\frac{0.99}{\omega_n}-C_n(C_*)^n\bigg)\frac{2}{\sqrt{2\pi}}\int_{\frac{1}{3}\tilde{C}_{n,p}}^\infty e^{-\frac{s^2}{2}} \ud s
$$
so that
$$
\mathbb P(\tau_b \leq \tau_g ~\text{or }t_{\tau_g} \geq \delta/3~\text{or }x_{\tau_g} \not \in A_{\delta/3}) \leq 1-\theta_{n,p}+\mathcal O(\epsilon)+\lambda.
$$
Denote $\theta:= \theta_{n,p}/2$ and thus, we have proven the claim for all $\epsilon$ and $\lambda$ small enough.
\end{proof}
\end{theorem}

If Player 1 plays by canceling the moves of the other player, we obtain Theorem \ref{paatulos p suurempi kaksi}. Observe that this strategy is not optimal for Player 1 in the sense that Player 1 also tries to cancel the moves that might benefit her.

The cancellation strategy was introduced in the paper \cite{luirops13} to prove Harnack's inequality for $p$-harmonic functions via tug-of-war games. In addition, the cancellation strategy can be used to prove regularity properties for viscosity solutions of the inhomogeneous $p$-Laplace equation (see \cite{ruosteenoja}).

\begin{theorem}\label{paatulos p suurempi kaksi}
If $y\in \partial \Omega$ satisfies the measure density condition, then it is game regular for $p>2$.
\begin{proof}
To establish the result, our aim is to use Lemma \ref{ehto} and therefore, to find a uniform lower bound for the probability that the game ends before exiting a given ball. If Player 1 plays by canceling the moves of the other player, the lower bound $\theta>0$ for the probability is obtained by using Theorem \ref{sylinteri arvio}. 

We can clearly assume that $y=0$. Let $\delta>0$, and consider the cylinder $B_{\delta/3}(0)\times [0,\delta/3]$. Define a constant $\delta_0$ as in \eqref{delta nolla}, and find $\epsilon_0>0$ and $\lambda>0$ small enough such that we can apply Theorem \ref{sylinteri arvio}. Let $x_0 \in B_{\delta_0}(0)$ and $\epsilon<\epsilon_0$. At every moment we can divide the game position as a sum of vectors
$$
x_0+\sum_{k\in I_1}v_k^1+\sum_{k\in I_2}v_k^2+\sum_{k\in I_3}v_k^3.
$$
Here, $I_1$ denotes the indices of rounds when Player 1 has moved with the vectors $v_k^1$ as her moves. Similarly, Player 2 has moved in the indices of rounds $I_2$ with the moves $v_k^2$ as his moves. The random movements have occured in the indices of rounds $I_3$, and these random movements are denoted by $v_k^3$. 

Let 
$$
M:=2\bigg \lceil \frac{ |x_0|}{\epsilon}\bigg \rceil, 
$$
where the factor 2 is due to the fact that the players cannot step to the boundary of $B_\epsilon (x_j)$ for any $j$. Define the following strategy $S_1^*$ for Player 1 for the game that starts from $x_0$. She always tries to cancel the earliest move of Player 2 which she has not yet been able to cancel. If all the moves at that moment are cancelled and she wins the coin toss, she moves the game point by the vector
$$
-\epsilon/2\frac{x_0}{|x_0|}.
$$
She does this until she has won $M-1$ more coin tosses than Player 2. If she wins her $M$th more coin toss, her move will be such that the game position is
$$
\sum_{k\in I_3}v_k^3
$$
after the move. Observe that the game, with the strategy $S_1^*$, is related to the cylinder walk, when we start the cylinder walk from the point $(0,M\epsilon/2)$ with $M\epsilon/2\to |x_0|<\delta_0$ as $\epsilon \to 0$.

Let us define three conditions for the game sequences of the game:

\begin{enumerate}[label=(\textbf{\Alph*})]
\item Player 1 has won the coin toss $M$ more times than Player 2, and at the moment this happens, the game sequence is in the set $\Omega^c$. 
\item  Player 2 has won the coin toss at least $ \frac{\delta}{3\epsilon}$ more times than Player 1. 
\item $|\sum_{k\in I_3}v_k^3|\geq \frac{\delta}{3}$.
\end{enumerate}
We are interested in the following event
$$
\textbf{X}:=\{\text{the condition $(\textbf{A})$ happens before conditions $(\textbf{B})$ and $(\textbf{C})$}\},
$$
and Theorem \ref{sylinteri arvio} states that there is a constant $\theta:=\theta_{n,p}>0$ such that
$$
\mathbb P_{S_1^*,S_2}(\textbf{X}) \geq \theta.
$$
Now, we can estimate
$$
\mathbb P_{S_1^*,S_2}(\text{the game ends before exiting the ball $B_\delta(0)$}) \geq \mathbb P_{S_1^*,S_2}(\textbf{X}).
$$
Above, we also used the fact that the game sequences for which the game has ended before Player 1 has won $M$ more coin tosses than Player 2 are good for our purposes. To finish the proof, we can use Lemma \ref{ehto}, and thus the proof is complete.
\end{proof}
\end{theorem}

It is worth mentioning that in the case $p>n$, every point becomes game regular. This is proved in \cite{peress08}, and the same also holds for the version of the game considered in this paper. Roughly, as $p$ increases, the probability for the player to end the game before exiting a given ball increases.

\appendix
\section{Hitting probabilities for a cylinder walk}

\renewcommand{\theequation}{\Alph{section}.\arabic{equation}}

Fix the cylinder size $r>0$. The cylinder walk in a cylinder $B_{r}(0)\times [0,r]\subset \mathbb R^{n+1}$ can be constructed by combining three independent random constructions. There is a 'horizontal' random walk with the initial position $\tilde{x}_0=x\in B_r(0)$. The point $\tilde{x}_{j+1}$ is chosen according to the uniform distribution in the ball $B_\epsilon(\tilde{x}_j)\subset \mathbb R^n$ for all $j\geq 0$. Further, there is a 'vertical' random walk in the real axis with steps $+\epsilon$ or $-\epsilon$ and with the initial position $\tilde{t}_0=t\in]0,r[$. The next positions are $\tilde{t}_{j+1}=\tilde{t}_{j}+\epsilon$ or $\tilde{t}_{j+1}=\tilde{t}_{j}-\epsilon$ both with probability $\frac{1}{2}$ for all $j\geq 0$. In addition, there is the increasing sequence
$$
U_j=\sum_{m=1}^j \ber_m,
$$
where the $\ber_m$:s are independent Bernoulli variables with $\ber_m(\omega) \in \{0,1\}$ and $\mathbb P(\ber_m=1)=\alpha\in]0,1[$. Thus, a copy of the cylinder walk is obtained by letting for $j\geq 0$
$$
t_j=\tilde{t}_{U_j},~~~x_j=\tilde{x}_{j-U_j}.
$$

Let $\tau_g$ stand for the first moment $t_j$ exits the cylinder through its bottom or top, and let $\tilde{\tau}_g$ stand for the first moment $\tilde{t}_j$ exits the cylinder through its bottom or top. 

Recall Hoeffding's (or Azuma's or Bernstein's) inequality for a sum of independent and identically distributed random variables (see for example \cite[p.\,198]{klenke08}).

\begin{theorem}\label{azuma}
Let $Y_m$ be independent and identically distributed symmetric $\mathbb R^n$-valued random variables, $m \in \{1,2,\dots, N\}$, that are uniformly bounded: $|Y_m|\leq b$ almost surely for all $m$. Then,
$$
\mathbb P \bigg(\max_{1\leq m\leq N} \Big|\sum_{i=1}^m Y_i\Big|\geq \lambda\bigg)\leq 4n\exp \bigg(- \frac{\lambda^2}{2Nb^2n}\bigg).
$$
\end{theorem}

In the theorem above, the factor 4 instead of 2 comes from the use of Levy-Kolmogorov's inequality (see for example \cite[p.\,397]{shiryaev96})
\begin{align*}
\mathbb P \bigg(\max_{1\leq m\leq N} \Big|\sum_{i=1}^m Y_i\Big|\geq \lambda\bigg) \leq 2 \mathbb P  \bigg( \Big|\sum_{i=1}^N Y_i\Big|\geq \lambda\bigg).
\end{align*}

We assume that $x=0$ i.e. $\tilde{x}_0=0$ and $\tilde{t}_0=t$ with $0<t<r$ and denote $\beta=1-\alpha$. 
\begin{lemma}\label{LEMMA A1}
Let $\tau_g$ and $\tilde{\tau}_g$ be as above. The random variables $\tilde{\tau}_g$ and $\tau_g-\tilde{\tau}_g$ have the following inequalities for all $a>0$
\begin{align}
\mathbb P(\tilde{\tau}_g \geq a\epsilon^{-1}) &\geq 1 -\mathcal O (\epsilon)\text{ and }\label{eka tulos}\\
\mathbb P(\tau_g-\tilde{\tau}_g \geq a\epsilon^{-1})& \geq 1 -\mathcal O (\epsilon)\label{toka tulos}.
\end{align}

\begin{proof}
The vertical movement consists of the moves $+\epsilon$ or $-\epsilon$ in the real axis. Let $Y_i$ be independent and identically distributed random variables  with $Y_i (\omega)\in\{-\epsilon,\epsilon\}$ and $\mathbb P(Y_i=\epsilon)=\mathbb P(Y_i=-\epsilon)=\frac{1}{2}$ for all $i$. Recall the cylinder size $B_{r}(0)\times [0,r]$. Now,
\begin{align*}
\mathbb P(\tilde{\tau}_g \geq a \epsilon^{-1})&=\mathbb P\bigg(\max_{k \leq a\epsilon^{-1}} \sum_{i=1}^k Y_i <\min \{t,r-t\} \bigg)\\
&=1- \mathbb P\bigg(\max_{k \leq a\epsilon^{-1}} \sum_{i=1}^k Y_i  \geq \min \{t,r-t\} \bigg).
\end{align*}
Random variables $Y_i$ are bounded, $|Y_i|\leq \epsilon$ for all $i \geq 1$. By using Hoeffding's inequality i.e. Theorem \ref{azuma}, we can deduce that
$$
\mathbb P\bigg(\max_{k \leq a\epsilon^{-1}} \sum_{i=1}^k Y_i  \geq \min \{t,r-t\} \bigg) \leq 4\exp\bigg( -\frac{(\min \{t,r-t\})^2}{2a\epsilon }\bigg) \leq \mathcal O(\epsilon).
$$
Consequently, we have proven the first part \eqref{eka tulos}.

For the second part, let us consider the event 
\begin{equation}\label{joukko B}
B:= \{0.99\alpha \tau_g<\tilde{\tau}_g<1.01\alpha\tau_g\}.
\end{equation}
Denote the sets
\begin{align*}
B^*&:=\{U_j<1.01\alpha j~\text{for all }j\geq a\epsilon^{-1}\}~\text{and} \\
B_*&:=\{U_j>0.99\alpha j~\text{for all }j\geq a\epsilon^{-1}\}.
\end{align*}
Again, apply Hoeffding's inequality with $Y_m=\ber_m-\alpha$, $\lambda= 0.01\alpha j$, $b=1$ and $N=j$ to get
\begin{align*}
\mathbb P(U_j \geq 1.01\alpha j)=\mathbb P(U_j-\alpha j \geq 0.01\alpha j) 
&\leq \mathbb P(|U_j-\alpha j| \geq 0.01\alpha j) \\
&\leq 4 \exp\bigg(-\frac{\alpha^2j}{2\cdot10^{4}}\bigg).
\end{align*} 
In a similar fashion, we can calculate 
$$
\mathbb P (U_j \leq 0.99\alpha j) \leq 4 \exp\bigg(-\frac{\alpha^2j}{2\cdot 10^{4}}\bigg).
$$
Thus, the summing over all indices $j\geq a\epsilon^{-1}$ leads to
\begin{align*}
\mathbb P\big((B^*)^c\big)&=\mathbb P (U_j \geq 1.01\alpha j~\text{for some }j\geq a\epsilon^{-1}) \\
&\leq \sum_{j\geq a\epsilon^{-1}}\mathbb P(U_j\geq 1.01\alpha j) \\
&\leq\sum_{j\geq a\epsilon^{-1}} 4 \exp\bigg(-\frac{\alpha^2j}{2\cdot 10^{4}}\bigg) \leq \mathcal O(\epsilon),
\end{align*}
and similarly, 
$$
\mathbb P\big(B_*\big) \geq 1-\mathcal O(\epsilon).
$$
These calculations imply 
$$
\mathbb P\big(B_*\text{ and } B^* \big)\geq \mathbb P(B_*)-\mathbb P\big((B^*)^c\big)\geq 1- \mathcal O(\epsilon).
$$
Observe that 
$$
\{B_*\text{ and } B^*\text{ and } \tau_g \geq a\epsilon^{-1}\} \subset B,
$$
and thus,
\begin{align*}
\mathbb P\big(B\big)&\geq \mathbb P(\tau_g \geq a\epsilon^{-1})-\mathcal O(\epsilon).
\end{align*}
Since $\tau_g \geq \tilde{\tau}_g$ always, we have $\{\tilde{\tau}_g \geq a\epsilon^{-1}\}\subset \{\tau_g \geq a\epsilon^{-1}\}$. Combining this with the first result  \eqref{eka tulos} we can deduce that
\begin{align}\label{alpha kerto}
\mathbb P(B)&\geq 1- \mathcal O(\epsilon).
\end{align}
In addition, we have by a direct calculation
$$
B \subset\bigg \{\frac{\beta-0.01\alpha}{1.01\alpha}\tilde{\tau}_g<\tau_g-\tilde{\tau}_g<\frac{\beta+0.01\alpha}{0.99\alpha}\tilde{\tau}_g \bigg \}.
$$
Therefore, we obtain by using the inequalities \eqref{eka tulos} and \eqref{alpha kerto} 
\begin{align*}
\mathbb P(\tau_g-\tilde{\tau}_g \geq a\epsilon^{-1})& \geq \mathbb P \bigg(\frac{\beta-0.01\alpha}{1.01\alpha}\tilde{\tau}_g \geq a\epsilon^{-1}\text{ and } B\bigg) \\
&=\mathbb P\big(\tilde{\tau}_g\geq 1.01\alpha a((\beta-0.01\alpha)\epsilon)^{-1}\text{ and }B\big) \\
&\geq \mathbb P\big(\tilde{\tau}_g\geq 1.01\alpha a((\beta-0.01\alpha)\epsilon)^{-1}\big)-\mathbb P\big(B^c\big) \\
&\geq1-\mathcal O(\epsilon).
\end{align*}
Thus, we have proven the second inequality \eqref{toka tulos}.
\end{proof}
\end{lemma}

\begin{lemma}\label{LEMMA A2}
Let $\tau_g, \tilde{\tau}_g,t,r$ and $\alpha,\beta >0$ such that $\alpha+\beta=1$ be as at the beginning of the appendix. Then, we have for all $a>0$ that
\begin{align*}
\mathbb P(\tau_g-\tilde{\tau}_g \geq a\epsilon^{-2}) \leq 1- \frac{2}{\sqrt{2 \pi}} \int_{\frac{\min\{t,r-t\}}{\sqrt{a}}\nu_{n,p}}^\infty e^{-\frac{s^2}{2}} \ud s+ \mathcal O(\epsilon)
\end{align*}
with the constant
$$
\nu_{n,p}:=2\sqrt{\frac{\beta+0.01\alpha}{0.99\alpha}}.
$$
\begin{proof}
By using the inequality \eqref{alpha kerto} and the inclusion after it, we can deduce 
\begin{align*}
\mathbb P(\tau_g-\tilde{\tau}_g \geq a\epsilon^{-2}) &\leq \mathbb P(\tau_g-\tilde{\tau}_g \geq a\epsilon^{-2} \text{ and } B)+ \mathbb P\big(B^c\big)\\
& \leq \mathbb P\bigg(\tilde{\tau}_g \geq \frac{0.99\alpha a}{(\beta+0.01\alpha)\epsilon^{2}}\bigg)+\mathcal O(\epsilon),
\end{align*}
where the set $B$ is defined as in \eqref{joukko B}.

We estimate the probability of the event $\{\tilde{\tau}_g \geq d\epsilon^{-2}\}$ for all $d>0$. Consider the following independent and identically distributed random variables: $Z_i(\omega) \in \{1,-1\}$, $\mathbb P (Z_i=-1)=\mathbb P (Z_i=1)=\half$ and $\mathbb E[Z_i]^2=1$ for all $i \geq 1$. For these random variables, we have the following equality (see for example \cite[p.\,351]{klenke08})
$$
\mathbb P \bigg(\max_{1\leq m\leq  N}\sum_{i=1}^m Z_i\geq l\bigg)=2\mathbb P  \bigg( \sum_{i=1}^{N}Z_i\geq l\bigg)-\mathbb P  \bigg( \sum_{i=1}^{N} Z_i= l\bigg)
$$
for all integers $N \geq 1$ and $l\geq 1$. Further, since $\mathbb E|Z_i|^3=1<\infty$ for all $i\geq 1$, we can use the Berry-Esseen theorem to determine the speed in the central limit theorem (see for example \cite[p.\,63]{shiryaev96}), and thus
\begin{align*}
 &2\mathbb P  \bigg( \sum_{i=1}^{N} Z_i\geq l\sqrt{N}\bigg)-\mathbb P  \bigg( \sum_{i=1}^{N} Z_i= l\sqrt{N}\bigg) & \\
&\geq \frac{2}{\sqrt{2 \pi}} \int_l^\infty e^{-\frac{s^2}{2}}~ds-\mathcal O(N^{-1/2}).
\end{align*}
Observe that for all $\epsilon>0$ small enough
$$
\frac{\big\lceil \min\{t,r-t\} \epsilon^{-1}\big\rceil}{\sqrt{\lfloor d \epsilon^{-2} \rfloor}} \leq \frac{2\min\{t,r-t\}}{\sqrt{d}}.
$$
Therefore, we have 
\begin{align*}
\mathbb P (\tilde{\tau}_g \geq d\epsilon^{-2})&\leq \mathbb P \bigg(\max_{1\leq m\leq \lfloor d\epsilon^{-2}\rfloor} \sum_{i=1}^m Z_i<\big\lceil \min\{t,r-t\} \epsilon^{-1}\big\rceil\bigg)\\
&\leq 1- \frac{2}{\sqrt{2 \pi}} \int_{\frac{2\min\{t,r-t\}}{\sqrt{d}}}^\infty e^{-\frac{s^2}{2}} \ud s+\mathcal O(\epsilon).\qedhere
\end{align*}

\end{proof}
\end{lemma}
Lemma \ref{LEMMA tod nak} is now an immediate consequence of Lemmas \ref{LEMMA A1} and \ref{LEMMA A2}.

Next, we prove a technical result (Lemma \ref{ehka kiinnostava} below) that we use in Section 3 above. First, in order to keep the calculations simple, let the dimension $n$ be one for now. Assume that $Z$ is distributed according to the uniform distribution in $]-\epsilon,\epsilon[$ for some $\epsilon>0$. Then for two independent $Z_1$ and $Z_2$ both distributed as $Z$, the density of the random variable $Z_1+Z_2$ can be computed via convolution. Thus, since $f_Z(x)=1/(2\epsilon)\1_{]-\epsilon,\epsilon[}(x)$, we have
$$
f_{Z_1+Z_2}(x)=\int_{-\infty}^\infty f_Z(x-y) f_Z(y) \ud y=\bigg(\frac{1}{2\epsilon}\bigg)^2\big(2\epsilon-|x|\big)\1_{]-2\epsilon,2\epsilon[}(x).
$$
For any $k\geq 1$, denote the density $f_k:=f_{\sum_{i=1}^k Z_i}$, where $Z_i$ are independent random variables distributed as $Z$. Similarly as in the case $k=2$, we can deduce and prove by induction (see for example \cite[p.\,197]{renyi70}) that for any $k \geq 1$
$$
f_k(x)=\begin{cases}
\frac{1}{(k-1)!(2\epsilon)^k}\sum_{j=0}^{\big\lfloor \frac{x+k\epsilon}{2\epsilon}\big\rfloor} (-1)^j\binom {k} {j}(x+k\epsilon-2j\epsilon)^{k-1},~&\text{if } x\in]-k\epsilon,k\epsilon[\\
0,~&\text{otherwise}.
\end{cases}
$$
Unfortunately, it is hard to get quantitative estimates from it.

There have been a lot of studies on the concentration function of a sum of independent random variables (see for example \cite{esseen68}). However, we are interested in the pointwise value of the function $f_k$ at the origin, and we will estimate the value by hand for the reader in a rather accessible way.

The characteristic function of the random variable $Z$ can be easily calculated,
\begin{align*}
\varphi_Z(t)=\frac{1}{2\epsilon}\int_{-\epsilon}^\epsilon e^{itx} \ud x=\frac{\sin (\epsilon t)}{\epsilon t}.
\end{align*}
Let $k\geq 2$. Because of the independence,
\begin{align*}
\varphi_{\sum_{i=1}^k Z_i}(t)=\bigg (\frac{\sin (\epsilon t)}{\epsilon t}\bigg )^k.
\end{align*}
Now, we have $\int_{-\infty}^\infty |\varphi_{\sum_{i=1}^k Z_i}(t)| \ud t <\infty$, so we can use the well-known inversion formula
$$
f_k(x)=\frac{1}{2 \pi} \int_{-\infty}^{\infty}e^{-itx}\varphi_{\sum_{i=1}^k Z_i}(t) \ud t.
$$
This inversion formula yields 
\begin{align*}
f_k(0)&= \frac{1}{ \pi}\int_{0}^{\epsilon^{-1}}\bigg (\frac{\sin (\epsilon t)}{\epsilon t}\bigg )^k \ud t+\frac{1}{\pi}\int_{\epsilon^{-1}}^{\infty} \bigg (\frac{\sin (\epsilon t)}{\epsilon t}\bigg )^k \ud t.
\end{align*}
Define
$$
h(z):=2\frac{1-\cos z}{z^2}
$$
so that we have for any $0 \leq m \leq 2\pi$
\begin{align}\label{arvio sini}
\frac{\sin z}{z} \leq 1-h(m)\frac{z^2}{6}
\end{align}
for all $|z|\leq m$. This inequality is true since the function $\sin z/z$ decreases for $0<z\leq\pi$ implying
$$
\bigg(\frac{\sin(m/2)}{m/2}\bigg)^2\leq \bigg(\frac{\sin (z/2)}{z/2}\bigg)^2
$$
for all $0<z\leq 2\pi$. This inequality yields
$$
1-\cos z-h(m)\frac{z^2}{2}\geq 0
$$
for all $0<z\leq 2\pi$ so we have the inequality \eqref{arvio sini}, since both sides of the inequality \eqref{arvio sini} are even functions. By using the inequality \eqref{arvio sini}, a change of variables formula and the inequality $1-z\leq e^{-z}$ for all $z\in \mathbb R$, we have
\begin{align*}
\frac{1}{ \pi}\int_{0}^{\epsilon^{-1}}\bigg (\frac{\sin (\epsilon t)}{\epsilon t}\bigg )^k \ud t & = \frac{1}{ \pi \epsilon}\int_{0}^{1}\bigg (\frac{\sin z}{z}\bigg )^k \ud z \\
& \leq\frac{1}{ \pi \epsilon}\int_{0}^{1}\bigg (1-h(1)\frac{z^2}{6}\bigg)^k \ud z \\
&\leq \frac{1}{ \pi \epsilon}\int_{0}^{1} e^{-\frac{z^2kh(1)}{6}} \ud z.
\end{align*}
Again, via changing the variables we derive
$$
\frac{1}{ \pi \epsilon}\int_{0}^{1} e^{-\frac{z^2kh(1)}{6}} \ud z \leq \frac{1}{\epsilon}\sqrt{\frac{6}{k \pi h(1)}}\frac{1}{\sqrt{2\pi} } \int_{0}^{\infty} e^{-\frac{u^2}{2}} \ud u= \sqrt{\frac{3}{ 2\pi h(1)}}\frac{1}{\sqrt{k}\epsilon}.
$$
Thus, we have estimated
$$
\frac{1}{ \pi}\int_{0}^{\epsilon^{-1}}\bigg (\frac{\sin (\epsilon t)}{\epsilon t}\bigg )^k\ud t \leq \sqrt{\frac{3}{ 2\pi h(1)}}\frac{1}{\sqrt{k}\epsilon}.
$$
Because $\sin z\leq 1$ for all $z\in \mathbb R$, we can estimate the second integral directly, and hence
\begin{align*}
\frac{1}{\pi}\int_{\epsilon^{-1}}^{\infty} \bigg (\frac{\sin (\epsilon t)}{\epsilon t}\bigg )^k\ud t& \leq \frac{1}{ \pi \epsilon}\int_{1}^{\infty} \frac{1}{z^k}\ud z=\frac{1}{\pi \epsilon (k-1)}.
\end{align*}
Therefore, we have derived the estimate
$$
f_k(0) \leq \sqrt{\frac{3}{2 \pi h(1)}}\frac{1}{\sqrt{k}\epsilon}+\frac{1}{\pi \epsilon (k-1)}.
$$

Next, we extend the argument to the higher dimensions as well. Assume that $Z$ is a random vector with the uniform distribution in the $n$-ball $B_\epsilon(0)$, $n\geq 1$. The density of the random vector $Z$ is 
$$
f_Z(x) = \frac{1}{|B_\epsilon(0)|} \1_{B_\epsilon(0)}(x).
$$
Using the same approach as in dimension one, we first need the characteristic function of the random vector $Z$. Denote the measure of the unit ball by $\omega_n:=|B_1(0)|=\pi^{n/2}/\Gamma(\frac{n}2+1)$, where the function $\Gamma$ is the usual gamma function. The random variable $Z$ is invariant under rotation i.e. the density function is a constant on every sphere $S^{n-1}_r(0):=\{x\in\mathbb R^n: |x|=r\}$ for all $r>0$. Hence, by rotating the ball $B_\epsilon(0)$, we see that $\varphi_Z(u)=\varphi_Z \big((r,0,\dots,0)\big)$ for all $u \in \mathbb R^n$ such that $|u|=r$. Let $r>0$, and direct computation with a change of variables $x=\epsilon y$ yields 
\begin{align*}
\varphi_Z \big((r,0,\dots,0)\big)&= \int_{\mathbb R^n} e^{ir x_1} f_Z(x) \ud x \\
& =\frac1{\om_n}\int_{B_1(0)} e^{i \eps r y_1 }  \ud y_1\ldots \ud y_n\\
&=\frac{\om_{n-1}}{\om_{n}}\int_{-1}^{1}(1-y^2_1)^{(n-1)/2} e^{i \eps ry_1 } \ud y_1 \\
&=\frac{\om_{n-1}}{\om_{n}}\int_{-1}^{1}(1-y^2_1)^{(n-1)/2} \cos( \eps ry_1 ) \ud y_1.
\end{align*}

A spherical Bessel function of order $n/2$, often denoted by $J_{n/2}(z)$, has an integral representation 
\begin{align*}\label{bessel integral}
J_{n/2}(z)=\bigg(\frac{z}{2}\bigg)^{n/2}\frac{1}{\Gamma(\frac{n+1}{2})\sqrt{\pi}}\int_{-1}^1(1-t^2)^{\frac{n-1}{2}}\cos(zt) \ud t
\end{align*}
(see for example \cite{watson44}). We can use this integral formula to obtain
\begin{align*}
\frac{\om_{n-1}}{\om_{n}}\int_{-1}^{1}(1-y^2_1)^{(n-1)/2} \cos( \eps ry_1 ) \ud y_1=(\eps r/2)^{-n/2}\Gamma(\frac{n}2+1)J_{n/2}(\eps r).
\end{align*}
Thus, we have derived the characteristic function
\begin{equation}\label{karakteristinen funktio}
\varphi_Z(u)=\bigg(\frac{2}{\eps |u|}\bigg)^{n/2}\Gamma(\frac{n}2+1)J_{n/2}\big(\eps |u|\big)
\end{equation}
for all $u\in \mathbb R^n$. Spherical Bessel functions have a connection to our calculations in dimension $n=1$, since one could show that
$$
\frac{\sin z}{z}=\sqrt{\frac{\pi}{2z}}J_{\frac{1}{2}}(z)
$$
holds for all $z\in \mathbb R$.

It is possible to express $J_{n/2}(z)$ as a product of factors such that each factor vanishes at one of the zeros of $z^{-n/2}J_{n/2}(z)$. Denote the zeros of the function $z^{-n/2}J_{n/2}(z)$ by $\pm j_{n/2,1},\pm j_{n/2,2}, \pm j_{n/2,3},\dots$ with $j_{n/2,l}>0$ for all $l =1,2,\dots$ and $j_{n/2,1}\leq j_{n/2,2}\leq j_{n/2,3}\leq \cdots$. Then, we have the infinite product formula of the Bessel function
\begin{equation}\label{bessel product}
J_{n/2}(z)=\bigg(\frac{z}{2}\bigg)^{n/2}\frac{1}{\Gamma(\frac{n}{2}+1)}\prod_{l=1}^\infty \bigg(1-\frac{z^2}{j_{n/2,l}^2}\bigg)
\end{equation}
(see \cite[p.\,497-498]{watson44}). The number of zeros of $z^{-n/2}J_{n/2}(z)$ between the origin and the point
$$
l_m:=m\pi+\frac{\pi}{4}\big(n+1\big)
$$
is exactly $m$ for all $m$ big enough (see \cite[p.\,495-497]{watson44}). Consequently, the infinite sum $\sum_{l=1}^\infty j_{n/2,l}^{-2} $ converges, since
$$
\sum_{l=p}^\infty j_{n/2,l}^{-2} \leq  \sum_{l=p}^\infty \bigg(\frac{1}{(l-1)\pi+\pi/4(n+1)}\bigg)^2<\infty
$$
for some $p$ big enough. Therefore, the infinite product in the formula \eqref{bessel product} is  well-defined for all $z\in \mathbb R$.

Via independence we have 
$$
\varphi_{\sum_{i=1}^k Z_i}(u) =\big(\varphi_Z (u)\big)^k, 
$$
and the inversion formula together with the characteristic function \eqref{karakteristinen funktio}, the infinite product formula \eqref{bessel product} and a change of variables $z=\epsilon u$ yield 
\begin{align*}
f_k(0)&=\frac{1}{(2\pi)^n} \int_{\mathbb R^n} \big(\varphi_{Z} (u)\big)^k\ud u \\
&=\frac{1}{(2\pi)^n\epsilon^n} \int_{B_s(0)}\bigg[\prod_{l=1}^\infty \bigg(1-\frac{ |z|^2}{j^2_{n/2,l}}\bigg)\bigg]^k\ud z\\
&+\frac{2^{kn/2}\Gamma(\frac{n}2+1)^k}{(2\pi)^n\epsilon^n} \int_{\mathbb R^n \setminus B_s(0)}\bigg(\frac{1}{|z|}\bigg)^{kn/2}\big(J_{n/2}\big( |z|\big)\big)^k \ud z
\end{align*}
for all $s>0$.

Now, the function 
$$
1-\frac{ |z|^2}{j^2_{n/2,l}}\geq 0
$$
for all $l\geq 1$, if $0\leq |z|\leq j_{n/2,1}$. In addition, since $1-z\leq e^{-z}$ for all $z\in \mathbb R$, we have
\begin{align*}
&\frac{1}{(2\pi)^n\epsilon^n} \int_{B_{j_{n/2,1}}(0)}\bigg[\prod_{l=1}^\infty \bigg(1-\frac{ |z|^2}{j^2_{n/2,l}}\bigg)\bigg]^k\ud z \\
&\leq  \frac{|S^{n-1}_1|}{(2\pi)^n\epsilon^n} \int_0^{j_{n/2,1}} e^{-r^2k\sum_{l=1}^\infty j^{-2}_{n/2,l}}r^{n-1}\ud r \\
&\leq \frac{|S^{n-1}_1|}{(2\pi)^n\epsilon^n} \int_0^\infty e^{-r^2k\sum_{l=1}^\infty j^{-2}_{n/2,l}}r^{n-1}\ud r. 
\end{align*}
Hence, we can integrate with a change of variables $r=\big(k\sum_{l=1}^\infty j^{-2}_{n/2,l}\big)^{-1/2}t$ to obtain
\begin{align*}
& \frac{|S^{n-1}_1|}{(2\pi)^n\epsilon^n} \int_0^\infty e^{-r^2k\sum_{l=1}^\infty j^{-2}_{n/2,l}}r^{n-1}\ud r \\
&=\frac{n \int_0^\infty e^{-t^2}t^{n-1}\ud t}{\Gamma(\frac{n}{2}+1)\pi^{n/2}2^{n}\big(\sum_{l=1}^\infty j^{-2}_{n/2,l}\big)^{n/2} }\bigg(\frac{1}{\sqrt{k}\epsilon}\bigg)^n.
\end{align*}
Thus, there is a constant 
$$
c_n^1:=\frac{n \int_0^\infty e^{-t^2}t^{n-1}\ud t}{\Gamma(\frac{n}{2}+1)\pi^{n/2}2^{n}\big(\sum_{l=1}^\infty j^{-2}_{n/2,l}\big)^{n/2} }>0
$$ 
such that
$$
\frac{1}{(2\pi)^n\epsilon^n} \int_{B_{j_{n/2,1}}(0)}\bigg[\prod_{l=1}^\infty \bigg(1-\frac{ |z|^2}{j^2_{n/2,l}}\bigg)\bigg]^k \ud z \leq c_n^1 \bigg(\frac{1}{\sqrt{k}\epsilon}\bigg)^n.
$$

The function $J_{n/2}\big(|z|\big)\leq 1$ for all $z\in \mathbb R^n$ (see for example \cite{watson44}). Therefore, we get by a direct calculus
\begin{align*}
&\frac{2^{kn/2}\Gamma(\frac{n}2+1)^k}{(2\pi)^n\epsilon^n} \int_{\mathbb R^n \setminus B_{j_{n/2,1}}(0)}\bigg(\frac{1}{|z|}\bigg)^{kn/2}\big(J_{n/2}\big( |z|\big)\big)^k \ud z \\
&\leq \frac{|S^{n-1}_1|2^{kn/2}\Gamma(\frac{n}2+1)^k}{(2\pi)^n\epsilon^n} \int_{j_{n/2,1}}^\infty r^{n-1-kn/2} \ud r \\
&=\frac{(j_{n/2,1})^n}{2^{n-1}\Gamma(\frac{n}2+1)\pi^{n/2}\epsilon^n}\bigg(\frac{1}{k-2}\bigg)\bigg(\bigg(\frac{2}{j_{n/2,1}}\bigg)^{n/2}\Gamma(\frac{n}2+1)\bigg)^k
\end{align*}
for all $k>2$. There exists the following lower bound for the first zero $j_{v,1}$ (see \cite{ifantiss85} and for example \cite{elbert01})
$$
j_{v,1}>v+\frac{\pi+1}{2}>v+2
$$
for all $v>- \half$. Thus, if $n$ is even, $n=2h$ for some $h\geq 1$, we get
$$
\frac{\Gamma(h+1)2^h}{(j_{h,1})^h}<\frac{h!2^h}{(h+2)^h}<1.
$$
Similarly, if $n$ is odd, $n=2h+1$ for some $h\geq 0$, we get
$$
\frac{\Gamma(h+\frac{3}{2})2^{h+1/2}}{(j_{h+1/2,1})^{h+1/2}}<\frac{(2h+2)!2^h\sqrt{2 \pi}}{4^{h+1}(h+1)!(h+2.5)^{h+1/2}}<1.
$$
Hence, there exists a constant $k_0:=k_{0,n}> 2$ such that
\begin{align*}
\bigg(\frac{1}{k-2}\bigg)\bigg(\bigg(\frac{2}{j_{n/2,1}}\bigg)^{n/2}\Gamma(\frac{n}2+1)\bigg)^k \leq \bigg(\frac{1}{\sqrt{k}} \bigg)^n
\end{align*}
for all $k\geq k_0$. Denote
$$
c_n^2:= \frac{(j_{n/2,1})^n}{2^{n-1}\Gamma(\frac{n}2+1)\pi^{n/2}}>0
$$
and
$$
C_n:=2\max \{c_n^1,c_n^2\}.
$$
Thus, we have derived the estimate
\begin{equation}\label{lopullinen arvio origo}
f_k(0)\leq C_n \bigg(\frac{1}{\sqrt{k}\epsilon} \bigg)^n
\end{equation}
for all $k\geq k_0$.

Let $k\geq k_0$. Theorem \ref{azuma} implies that there is a constant $C_1:=C_{1,n}>0$ big enough such that for all $\epsilon>0$
$$
\mathbb P\bigg(\Big|\sum_{i=1}^k Z_i\Big| < C_1\sqrt{k}\epsilon \bigg) \geq 0.99.
$$
By using the convolution formula, we have that
\begin{align}
\begin{split}\label{konvoluution avulla}
f_k(x)=\int_{\mathbb R^n}f_{k-1}(x-y)\1_{B_\epsilon(0)}(y)\ud y&=\int_{B_\epsilon(0)}f_{k-1}(x-y)\ud y\\
&=\int_{B_\epsilon(x)}f_{k-1}(y)\ud y
\end{split}
\end{align}
holds for all $x\in \mathbb R^n$. The function $f_1$ is a decreasing radial function. Thus, we can deduce by using the formula \eqref{konvoluution avulla} that $f_2$ is also a decreasing radial function, and by induction $f_k$ as well. Therefore, we can denote the density $f_k$ as a function of the radius $|u|$ for all $u \in \mathbb R^n$, and we have for any $C_*\in ]0,C_1[$ 
$$
f_k(0)|B_{C_*\sqrt{k}\epsilon}(0)|+f_k\big(C_*\sqrt{k}\epsilon\big)\Big(|B_{C_1\sqrt{k}\epsilon}(0)|-|B_{C_*\sqrt{k}\epsilon}(0)|\Big)\geq 0.99.
$$
This inequality yields
\begin{align*}
f_k\big(C_*\sqrt{k}\epsilon\big) &\geq \frac{0.99-f_k(0)|B_{C_*\sqrt{k}\epsilon}(0)|}{|B_{C_1\sqrt{k}\epsilon}(0)|-|B_{C_*\sqrt{k}\epsilon}(0)|} \\
&\geq \frac{0.99}{|B_{C_1\sqrt{k}\epsilon}(0)|}-f_k(0)\bigg(\frac{C_*}{C_1}\bigg)^n.
\end{align*}
Now, we use the estimate \eqref{lopullinen arvio origo} to obtain
\begin{align*}
f_k\big(C_*\sqrt{k}\epsilon\big)&\geq \frac{0.99}{|B_{C_1\sqrt{k}\epsilon}(0)|}-C_n \bigg(\frac{C_*}{C_1 \sqrt{k}\epsilon}\bigg)^n \\
&=\bigg(\frac{1}{C_1}\bigg)^n\bigg(\frac{0.99}{\omega_n}-C_n(C_*)^n\bigg)\bigg(\frac{1}{\sqrt{k}\epsilon}\bigg)^n.
\end{align*}
Thus, we have proven the following lemma.
\begin{lemma}\label{ehka kiinnostava}
Let $\epsilon>0$ and let $Z$ be distributed according to the uniform distribution in the ball $B_\epsilon(0)\subset \mathbb R^n$. For any $k\geq 2$, denote the density of the random variable $\sum_{i=1}^k Z_i$ by $f_k$, where the random variables $Z_i$, $i\in \{1,\dots,k\}$, are independent and distributed as $Z$. Then $f_k$ is a decreasing radial function, and there exist universal constants $k_0:=k_{0,n}>2$, $C_1:=C_{1,n}>0$ and $C_n>0$ such that for all $k \geq k_0$ and $C_*\in ]0,C_1[$ we have
\begin{equation}\label{lopullinen}
f_k\big(C_*\sqrt{k}\epsilon\big)\geq \bigg(\frac{1}{C_1}\bigg)^n\bigg(\frac{0.99}{\omega_n}-C_n(C_*)^n\bigg)\bigg(\frac{1}{\sqrt{k}\epsilon}\bigg)^n.
\end{equation}
\end{lemma}
Observe that 
$$
f_k\big(C_*\sqrt{k}\epsilon\big) \geq \zeta\bigg(\frac{1}{\sqrt{k}\epsilon}\bigg)^n
$$ 
for some $\zeta:=\zeta_n>0$, if we choose $C_*>0$ so small that
$$
C_*<\bigg(\frac{0.99}{\omega_nC_n}\bigg)^{1/n}.
$$


\def\cprime{$'$} \def\cprime{$'$}

\end{document}